\pdfoutput=1
\RequirePackage{ifpdf}
\ifpdf 
\documentclass[pdftex]{sigma}
\else
\documentclass{sigma}
\fi

\begin{document}


\renewcommand{\thefootnote}{$\star$}

\renewcommand{\PaperNumber}{076}

\FirstPageHeading

\ShortArticleName{Recursion Operators and Frobenius Manifolds}

\ArticleName{Recursion Operators and Frobenius Manifolds\footnote{This
paper is a contribution to the Special Issue ``Geometrical Methods in Mathematical Physics''. The full collection is available at \href{http://www.emis.de/journals/SIGMA/GMMP2012.html}{http://www.emis.de/journals/SIGMA/GMMP2012.html}}}

\Author{Franco MAGRI}

\AuthorNameForHeading{F.~Magri}

\Address{Dipartimento di Matematica ed Applicazioni, Universit\`{a} degli Studi di di Milano Bicocca,\\ Via Roberto Cozzi 53, 20125 Milano, Italy}
\Email{\href{mailto:franco.magri@unimib.it}{franco.magri@unimib.it}}

\ArticleDates{Received June 01, 2012, in f\/inal form October 05, 2012; Published online October 19, 2012}

\Abstract{In this note I exhibit a ``discrete homotopy'' which joins the category of F-manifolds to the category of Poisson--Nijenhuis manifolds, passing through the category of Frobenius manifolds.}


\Classification{35D45; 53D17; 37K10}

\renewcommand{\thefootnote}{\arabic{footnote}}
\setcounter{footnote}{0}

 \section*{Introduction}
In this short note I wish to present some ideas of a recent comparative study of three dif\/ferent types of manifolds  which are of interest in dif\/ferent areas of physics and mathematics. They are:
\begin{enumerate}\itemsep=0pt
\item[1)] F-manifolds;
\item[2)] Frobenius manifolds;
\item[3)] bi-Hamiltonian manifolds;
\end{enumerate}
The study is addressed to emphasize a point of view allowing to see these three dif\/ferent kinds of manifolds  under a unif\/ied perspective, almost as three dif\/ferent instances (or realizations) of a single geometrical idea.

Let us start by recalling f\/irst the def\/initions of the manifolds  we are interested in.

\textbf{F-manifolds.} F-manifolds  have been introduced by Hertling and Manin in 1999~\cite{1}. They are manifolds  where it is allowed to multiply vector f\/ields according to a rule of the form:
\[
\partial_{j} \circ \partial_{k} = C_{jk}^{l} \partial_{l}.
\]
The multiplication is assumed to be commutative, associative and with unity. The structure constant~$C_{jk}^{l}$ are the components of a symmetric vector-valued 2-form $C:TM \times TM \to TM,$ and the unity is a special vector f\/ield $e:M \to TM$. According to Hertling and Manin the manifold $M$ is a F-manifold if the multiplication $C$ and the unity $e$ verify the conditions~\cite{2}:
\begin{enumerate}\itemsep=0pt
\item[1)] $\mathrm{Lie}_{e}(C) = 0$;

\item[2)] $\mathrm{Lie}_{X \circ Y}(C) = X \circ \mathrm{Lie}_Y (C) + \mathrm{Lie}_X (C) \circ Y$.
\end{enumerate}

\textbf{Frobenius manifolds.}  Frobenius manifolds are manifolds  where each tangent space has the structure of a Frobenius algebra. The notion of Frobenius manifold has been introduced by Boris Dubrovin in 1992 \cite{3}. According to him, a manifold equipped with a multiplication $C:TM \times TM \to TM$, with a unity $e:M \to TM$, and with a metric $g: TM \times TM \to \mathbb{R}$ (not necessarily positive def\/inite) is a Frobenius manifolds if:
\begin{enumerate}\itemsep=0pt
  \item[1)] the metric is f\/lat: $\mathrm{Riemann}(g)=0$;
  \item[2)] the unity is covariantly constant: $\nabla e = 0$;
  \item[3)] the multiplication is symmetric with respect to the metric: $g(X, Y \circ Z) = g(X \circ Y, Z)$;
  \item[4)] the covariant dif\/ferential of the symmetric covariant tensor f\/ield $c(X,Y,Z) = g(X, Y \circ Z)$ is symmetric, that is it verif\/ies the symmetry condition: $\nabla_X c(Y, Z, W) = \nabla_Y c(X,Z,W)$.
\end{enumerate}
It can be shown that the axioms of Dubrovin entail that the multiplication $C$ and the unity $e$ verify the axioms of Hertling and Manin. Therefore, Frobenius manifolds are a particular class of F-manifolds.

\textbf{Bi-Hamiltonian manifols.}  Bi-Hamiltonian manifolds have been introduced in 1979. They are manifolds  equipped with a pair of compatible Poisson brackets. This means that the pencil
\[
\{ f,g \}_{\lambda} = \{ f,g \}_{1} + \lambda \{ f,g \}_{2}, \qquad \lambda \in \mathbb{R},
\]
def\/ined by the pair of Poisson brackets $\{ f,g \}_{1}$, and $\{ f,g \}_{2}$ satisf\/ies the Jacobi identity for any value of the parameter~$\lambda$. If one of the Poisson bracket is symplectic, the bi-Hamiltonian manifold is called a Poisson--Nijenhuis manifold. This type of manifolds  has been thoroughly studied in~\cite{4}.

At f\/irst sight, it is dif\/f\/icult to see any point of contact between the def\/initions of Frobenius manifolds and of Poisson--Nijenhuis manifolds. They belong to dif\/ferent geometries. Frobenius manifolds belong to Riemannian geometry, while Poisson--Nijenhuis manifolds  belong to symplectic geometry. Nevertheless it is known, after the work of Witten and Kontsevich, that the physical theories which are behind them, namely the topological f\/ield theories and the theory of integrable systems, are closely related. This occurrence leads to suspect that also between the correlated geometrical models there should exist strong points of contact. This idea has motivated the present research.

To attain this unif\/ication, I follow a technique of homotopic deformation. I imagine the three dif\/ferent categories of manifolds  as three distinct points in an abstract space (the space of all categories of smooth manifolds), and I conceive the idea to search for a simple path connecting these points. The construction of a path which joins the category of F-manifolds  to the category of Poisson--Nijenhuis manifolds, passing trough the category of Frobenius manifolds is the main novelty contained in the present paper.

\begin{figure}[h!]
\centering
\includegraphics{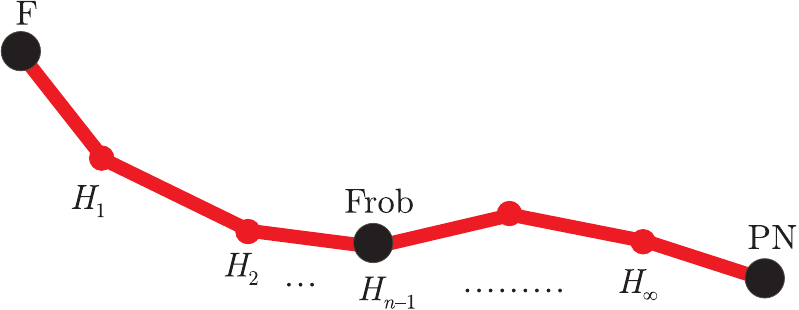}
\end{figure}

As shown by the f\/igure, this path is actually a polygonal line, whose vertices are denoted by $H_{1}, H_{2}, \dots, H_{\infty}$. Each vertex is a new category of manifolds, def\/ined according to a regular rule. Moving to the right, towards PN-manifolds, the complexity of the manifolds  increases. The manifolds  of type $H_{m+1}$ have indeed a richer geometrical structure than the manifolds  of type $H_{m}$, in the sense that they are manifolds  of type $H_{m}$ endowed with additional geometric proprieties. Therefore the vertices of the polygonal line represent a chain of categories of smooth manifolds obeying the inclusion relations:
\[
H_{1} \supset H_{2} \supset \cdots \supset H_m \supset H_{m+1} \supset \cdots \supset H_{\infty}.
\]
F-manifolds  and PN-manifolds  are the terminal points of this chain, since F-manifolds  are manifolds of type $H_1$ deprived of one structure, while PN-manifolds  are manifolds  of type $H_{\infty}$ enriched with an additional structure. They are not included into the chain, because their inclusion would destroy the validity of the simple recursive procedure allowing to pass from the category $H_m$  to the category $H_{m+1}$. To keep the simplicity of the scheme, they must be considered as limit points not belonging to the chain.

This paper has two specif\/ic aims:
\begin{enumerate}\itemsep=0pt
\item[1)] to explain the def\/inition of the new categories of manifolds  $H_m$;
\item[2)] to show that our polygonal line passes through the category of Frobenius manifolds, by a~concrete example.
\end{enumerate}

\section*{Manifolds of type $\boldsymbol{H_1}$}
The def\/inition of manifolds  of type $H_1$ comes directly from the theory of bi-Hamiltonian manifolds.

Let $M$ be a bi-Hamiltonian manifold, and let $X$ be a bi-Hamiltonian vector f\/ield def\/ined over it. If $P_1 : T^{*} M \to TM$ and $P_2 : T^{*} M \to TM$
are the Poisson bivectors associated with the pair of Poisson brackets giving $M$ the structure of bi-Hamiltonian manifold, the last requirement means that there exist on $M$ a pair of functions $h_1 : M \to \mathbb{R}$ and $h_2 : M \to \mathbb{R}$ such that
\[
X = P_1 \mathrm{d}h_1 = P_2 \mathrm{d}h_{2}.
\]
Assuming that $P_1$ is invertible, let us set
\[
K = P_2 P_{1}^{-1},
\]
and let us consider the manifold $M$ equipped with the recursion operator $K$, the bi-Hamiltonian vector f\/ield~$X$, and with the privileged 1-form $\theta = \mathrm{d}h_{1}$. It is well-known, within the theory of bi-Hamiltonian manifolds, that these three geometric objects satisfy the conditions:
\begin{enumerate}\itemsep=0pt
\item[1)] $\mathrm{Torsion}(K) = 0$;
\item[2)] $\mathrm{d} \theta = 0$;
\item[3)] $\mathrm{d}(K \theta) = 0$;
\item[4)] $\mathrm{Lie}_X (K)$ = 0;
\item[5)] $\mathrm{Lie}_X (\theta) = 0$.
\end{enumerate}

To pass from PN-manifolds  to the manifolds  of type $H_1$, one has to modify these conditions in two ways. First, one has to replace the strong condition
\[
\mathrm{Torsion}(K) = 0
\]
by the weaker condition
\[
\mathrm{Haantjes}(K) = 0.
\]
Secondly, one has to balance the weakening of the condition on $K$ by an additional condition on~$\theta$. This condition requires that the 1-form $\theta$ annihilates the torsion of $K$, which is now dif\/ferent from zero. So, the new condition is
\[
\theta (\mathrm{Torsion}(K)) = 0.
\]

\begin{definition}
A manifold $M$ equipped with a privileged vector f\/ield $X: M \to TM$, with a~privileged 1-form $\theta : M \to T^{*}M$, and with a~recursion operator $K: TM \to TM$ is a manifold of type $H_1$ if:
\begin{enumerate}\itemsep=0pt
\item[1)] $\mathrm{Haantjes}(K) = 0$;
\item[2)] $\mathrm{d}\theta = 0$;
\item[3)] $\mathrm{d} (K \theta) = 0$;
\item[4)] $\theta (\mathrm{Torsion}(K)) = 0$;
\item[5)] $\mathrm{Lie}_{X}(K) = 0$;
\item[6)] $\mathrm{Lie} _{X} (\theta) = 0$.
\end{enumerate}
\end{definition}

The additional condition, which demands that $\theta$ annihilates the torsion of~$K$, has a clear geometrical meaning. Together with the conditions $\mathrm{d}\theta = 0$ and $\mathrm{d} (K \theta) = 0$, it entails that the second iterated 1-form $K^2 \theta$ is closed as well. Thus, on any manifold of type $H_1$ there is a triple of 1-forms $(\theta, K \theta, K^2 \theta)$ which are closed, and therefore locally exact. This triple will be referred to as the \emph{short Lenard chain} associated with the recursion operator $K$ and the privileged 1-form $\theta$. By denoting by $K_0 : TM \to TM$ the identity map on the tangent bundle, and by renaming~$K_1$ the recursion operator $K$, one may condense the three conditions def\/ining a short Lenard chain in the more compact and neat form
\[
\mathrm{d}(K_{j} K_{l} \theta) = 0, \qquad j, l = 0, 1.
\]
This equation will be referred to as the def\/inition of the short Lenard chain of a manifolds  of type~$H_1$. Making use of this equation, one may give the following alternative def\/inition of manifolds  of type~$H_1$.

\begin{definition}
The manifold $M$ is a manifold of type $H_1$, if
\begin{enumerate}\itemsep=0pt
\item[1)] $\mathrm{Haantjes}(K_j) = 0$;
\item[2)] $\mathrm{d}(K_{j} K_{l} \theta) = 0$;
\item[3)] $\mathrm{Lie}_{X}(K_{j}) = 0$;
\item[4)] $\mathrm{Lie}_{X} (\theta) = 0$
\end{enumerate}
for $j, l = 0,1 $, where $ X $ is a privileged vector f\/ield and $ K_0 $ the identity map.
\end{definition}

To see the relation between manifolds  of type $H_1$ and F-manifolds, it is enough to concentrate the attention on the standard Lenard chain generated by the recursion operator $K$ and by the privileged vector f\/ield $X$. This chain is composed by the vector f\/ields
\[
X_0 = X, \quad X_1 = K X, \quad X_2 = K^2 X, \quad \dots
\]
Let us assume that the f\/irst $n=\mathrm{dim} M$ vector f\/ields of this chain are linearly independent, at least in some open region of $M$. Under the semisimplicity assumption
\[
X_0 \wedge X_1 \wedge \cdots \wedge X_{n-1} \neq 0,
\]
the Lenard chain def\/ines a frame $\{ X_{j} \}_{j=0,\dots , n-1}$ and a coframe $\{ \varepsilon^{j} \}_{j=0,\dots , n-1}$ on $M$. By using the coframe $\varepsilon^{j}$ and the powers $K^{l}$ of the recursion operators, one may then def\/ine the third-order tensor f\/ield
\[
C = \varepsilon^{0} \otimes Id + \varepsilon^{1} \otimes K + \varepsilon^{2} \otimes K^2 + \cdots + \varepsilon^{n-1} \otimes K^{n-1}.
\]
It def\/ines a multiplicative structure on the tangent bundle. It can be shown that the tensor f\/ield~$C$, and the privileged vector f\/ield~$X$ obey the axioms of Hertling and Manin. Therefore each semisimple manifold of type~$H_1$ is a semisimple F-manifold of Hertling and Manin. The converse is also true: any semisimple F-manifold is a semisimple manifold of type~$H_1$, stripped of the privileged 1-form~$\theta$.

\section*{Manifolds of type $\boldsymbol{H_2}$ and beyond}

Having replaced the condition $\mathrm{Torsion}(K) = 0$ by the weaker condition $\mathrm{Haantjes}(K) = 0$, one looses an important property of bi-Hamiltonian manifolds. On these manifolds  any short Lenard chain may be indef\/initely prolonged into an inf\/inite chain of closed 1-forms, to be referred to as the long Lenard chain. The process of extension is rather simple: it exploits the successive powers of the recursion operators. It is well-known, indeed, that on a bi-Hamiltonian manifold all the 1-form $K^{j} \theta$ are closed. So, there is not end to a Lenard chain on a bi-Hamiltonian manifold.

This property fails on a manifold of type~$H_{1}$. Normally it is not possible to prolong a short Lenard chain beyond its normal length by the iterated action of the recursion operator~$K$. However, this dif\/f\/iculty may be circumvented. To construct short Lenard chains composed of more than three closed 1-form, the strategy is to replace the powers of the recursion operator by new recursion operators suitably selected. On imitation of the powers of~$K$, the new recursion operators are assumed to commute. In its simplest form, this idea leads to the following def\/inition of manifolds  of type~$H_2$.

\begin{definition}
A manifold $M$ equipped with a privileged vector f\/ield $X$, a privileged 1-form $\theta$, and two recursion operators~$K_1$ and~$K_2$ is a manifold of type~$H_2$ if:
\begin{enumerate}\itemsep=0pt
\item[1)] $\mathrm{Haantjes}(K_j) = 0$;
\item[2)] $[K_j , K_l] = 0$;
\item[3)] $\mathrm{d}(K_j K_l \theta) = 0$;
\item[4)] $\mathrm{Lie}_{X}(K_{j}) = 0$;
\item[5)] $\mathrm{Lie}_{X} (\theta) = 0$
\end{enumerate}
for any value of the indices $j, l = 0, 1, 2$.  As usual, $K_0$ is the identity map on the tangent bundle~$TM$.
\end{definition}

The def\/inition of manifold of type $H_m$ in presently completely obvious. Each time one adds a new recursion operator, and one insists that the above conditions be satisf\/ied for any possible choice of the values of the indices $j$, $l$, ranging between~$0$ and~$m$.

\section*{WDVV equations}

Frobenius manifolds belong to the polygonal line def\/ined by the vertices $H_{1}, H_{2}, \dots, H_{\infty}$. Indeed, one may prove that a Frobenius manifold of dimension $n$ is a manifold of type $H_{n-1}$. Conversely, under very mild conditions of nondegeneracy, one may also prove that a manifold of type $H_{n-1}$ and dimension $n$ is a Frobenius manifold. The metric $g$ and the multiplication $C$ are reconstructed by using only the recursion operators $(K_{0}, K_{1}, K_{2}, \dots, K_{n-1})$, the privileged vector f\/ield~$X$, and the privileged 1-form $\theta$ of the manifold of type $H_{n-1}$. Thus the recursion operators have in the theory of Frobenius manifolds the same role they have in the theory of Poisson--Nijenhuis manifolds.

To prove these claims is outside the limits of the present note. For this reason, I limit myself to consider a concrete example. In this section I will study manifolds  of type~$H_{2}$ and dimension~$3$, and I will show that they are directly related to the WDVV equations, so important within the theory of Frobenius manifolds.

Let us start with the remark that a short Lenard chain on a manifold of type $H_2$ is composed by six closed 1-forms. Locally I can consider these forms as exact, setting:
\begin{alignat*}{4}
& \theta = \mathrm{d}A, \qquad &&  K_1 \theta = \mathrm{d}B, \qquad && K_2 \theta = \mathrm{d}C, &\\
& K_{1}^{2} \theta = \mathrm{d}P, \qquad &&  K_1 K_2 \theta = \mathrm{d}Q , \qquad && K_{2}^{2} \theta = \mathrm{d}R. &
\end{alignat*}
Let us  assume that the f\/irst three 1-form $\mathrm{d}A$, $\mathrm{d}B$, $\mathrm{d}C$ are linearly independent. The corres\-pon\-ding potentials are, therefore, local coordinates on the manifolds  $M$. Let us evaluate the action of the recursion operators $K_1$ and $K_2$ on this coordinate basis. One obtains
\begin{alignat*}{3}
& K_1 \mathrm{d}A = \mathrm{d}B, \qquad && K_2 \mathrm{d}A = \mathrm{d}C, & \\
& K_1 \mathrm{d}B = \mathrm{d}P, \qquad &&  K_2 \mathrm{d}B = \mathrm{d}Q, & \\
& K_1 \mathrm{d}C = \mathrm{d}Q,   \qquad && K_2 \mathrm{d}C = \mathrm{d}R, &
\end{alignat*}
where $P$, $Q$, $R$ are, so far, arbitrary functions of the coordinates $A$, $B$, $C$. These formulas def\/ine the most general pair of recursion operators generating a short Lenard chain of 1-forms on a~manifolds  of type~$H_2$ and dimension~$3$.

Let us now require that $K_1$ and $K_2$ commute. This condition restricts the choice of the arbitrary functions $P$, $Q$, $R$. The commutativity conditions are:
\[
K_2 \mathrm{d}P = K_1 \mathrm{d}Q, \qquad
 K_2 \mathrm{d}Q = K_1 \mathrm{d}R.
 \]
They imply that the functions $P$, $Q$, $R$ must solve the system of partial dif\/ferential equations
\begin{gather*}
\dfrac{\partial P}{\partial A}   = \dfrac{\partial P}{\partial C} \left( \dfrac{\partial Q}{\partial B} - \dfrac{\partial R}{\partial C} \right) + \dfrac{\partial Q}{\partial C} \left( \dfrac{\partial Q}{\partial C} - \dfrac{\partial P}{\partial B} \right), \\
\dfrac{\partial Q}{\partial A}   = \dfrac{\partial P}{\partial C} \dfrac{\partial R}{\partial B} - \dfrac{\partial Q}{\partial C} \dfrac{\partial Q}{\partial B}, \\
\dfrac{\partial R}{\partial A}   = \dfrac{\partial Q}{\partial B} \left( \dfrac{\partial Q}{\partial B} - \dfrac{\partial R}{\partial C} \right) + \dfrac{\partial R}{\partial B} \left( \dfrac{\partial Q}{\partial C} - \dfrac{\partial P}{\partial B} \right).
\end{gather*}

Another restriction comes from the existence of the vector f\/ield $X$. The symmetry conditions
\[
\mathrm{Lie}_X (K_j) = 0, \qquad  \mathrm{Lie}_X (\theta) = 0
\]
entail that the functions $A$, $B$, $C$, $P$, $Q$, $R$ have constant derivatives along~$X$. In particular, this means that the vector f\/ield $X$ has constant components on the frame def\/ined by the coordinates $A$, $B$, $C$. Thus
\[
X = a \dfrac{\partial}{\partial A} + b \dfrac{\partial}{\partial B} + c \dfrac{\partial}{\partial C},
\]
with $a$, $b$, $c$ constant parameters. Let us choose, for the sake of simplicity, $a = 0$, $b = 0$, $c = 1$, so that
\[
X = \dfrac{\partial}{\partial C}.
\]
The symmetry conditions then entail that the function $P$, $Q$, $R$ have the form
\[
P = \lambda C + \phi (A,B), \qquad Q = \mu C + \psi (A,B), \qquad R = \nu C + \chi (A,B),
\]
where $\lambda$, $\mu$, $\nu$ are arbitrary constants and the coordinate $C$ has been separated. Still for simplicity, let us choose $\lambda = 1$, $\mu = 0$, $\nu = 0$. By inserting
\[
P = C + \phi (A,B), \qquad Q =  \psi (A,B), \qquad R =    \chi (A,B)
\]
into the above partial dif\/ferential equations, one is left with
\[
\phi_{A} = \psi_{B}, \qquad \psi_{A} = \chi_{B}, \qquad \chi_{A} = \psi_{B}^{2} - \chi_B \phi_B.
\]
The solution is
\[
\phi = F_{AA}, \qquad \psi = F_{AB}, \qquad \chi = F_{BB}
\]
with
\[
F_{AAA} + F_{AAB}F_{BBB} - F_{ABB}^2 = 0.
\]

This is one of the classical WDVV equations for Frobenius manifolds of dimension~3~\cite{3}.
This example shows that this WDVV equation comes out naturally from the study of short Lenard chains on manifolds  of type~$H_2$ and dimension~3. The general WDVV equations come analogously from the study of Lenard chain on manifolds  of type $H_{n-1}$ and dimension~$n$.

\subsection*{Acknowledgements}

This note is a f\/irst account of a lasting research work done in collaboration with B.~Konopel\-chenko. To him I address my warm acknowledgements for the permission to use part of the common work to prepare the conference and this note. Many thanks are also due to B.~Dubrovin for the kind invitation to attend the conference on Geometrical Methods in Mathematical Physics.

\pdfbookmark[1]{References}{ref}
\LastPageEnding

\end{document}